\documentclass[12pt]{amsart}
\usepackage[final]{epsfig}
\usepackage{graphics}
\usepackage{amsmath}
\usepackage{amsfonts}
\usepackage{latexsym}
\usepackage{amssymb}
\usepackage{epstopdf}
\usepackage{graphicx}
\usepackage{url}
\usepackage[margin=1in]{geometry}

\newtheorem{lemma}{Lemma}
\newtheorem{proposition}[lemma]{Proposition}

\newtheorem{example}[lemma]{Example}
\newtheorem{theorem}{Theorem}
\newtheorem{definition}[lemma]{Definition}
\newtheorem{corollary}[lemma]{Corollary}

\newtheorem{conjecture}[theorem]{Conjecture}

\newcommand{\N}{{\mathbb N}}
\newcommand{\Z}{{\mathbb Z}}

\def\proof{\paragraph{Proof.}}
\newcommand{\proofend}{$\Box$\bigskip} 

\title{Monotonicity of the optimal perimeter in isoperimetric problems on ${\mathbb{Z}}^{k} \times {\mathbb{N}}^{d}$}

\begin{document}

\author{Emmanuel Tsukerman}
\address{Stanford University, Stanford, CA 94305}
\email{emantsuk@stanford.edu}

\date{}

\begin{abstract}
We prove general theorems for isoperimetric problems on lattices of the form ${\mathbb{Z}}^{k} \times {\mathbb{N}}^{d}$  which state that the perimeter of the optimal set is a monotonically increasing function of the volume under certain natural assumptions, such as local symmetry or being induced by an $\ell_p$-norm. The proved monotonicity property is surprising considering that solutions are not always nested (and consequently standard techniques such as compressions do not apply). The monotonicity results of this note apply in particular to vertex- and edge-isoperimetric problems in the $\ell_p$ distances and can be used as a tool to elucidate properties of optimal sets. As an application, we consider the edge-isoperimetric inequality on the graph ${\N}^2$ in the $\ell_\infty$-distance. We show that there exist arbitrarily long consecutive values of the volume for
which the minimum boundary is the same.
\end{abstract}

\maketitle

\section{Background}

Isoperimetric problems are classical objects of study in mathematics.
In general, such problems ask for sets whose boundary is smallest
for a given volume. A classic example dating back to ancient Greece
is to determine the shape in the plane for which the perimeter
is minimized subject to a volume constraint. 

Two notions of the boundary which are commonly considered are the vertex boundary
and the edge boundary. These are defined as follows. Let $A$
be a set of vertices in a graph $G=(V,E)$. The vertex boundary of
$A$ is defined to be the set $\{v\in V\setminus A:v\sim a\mbox{ for some }a\in A\}$.
In contrast, the \textit{edge boundary} of $A$ is defined to be the
set $\{e\in E:e=\{v,a\}\mbox{ for some }a\in A\mbox{ and }v\in V\setminus A\}$.
The vertex (edge) isoperimetric problem\textit{ }for $G$ asks for the minimum possible cardinality of the
vertex (edge) boundary for a $k$-element subset of $V$ for each
$k\in\mathbb{N}$. 

One of the settings in which discrete vertex- and edge-isoperimetric problems
have previously been studied is the graph of the infinite grid viewed in
the $\ell_{1}$ metric, in which two nodes are adjacent whenever their
distance in the $\ell_{1}$ metric is $1$ (see \cite{WangWang}
for the vertex-isoperimetric problem and \cite{BollobasLeader}
for the edge-isoperimetric problem). Natural analogues of these problems
are the vertex- and edge-isoperimetric problems on the infinite grid
in the $\ell_{\infty}$ metric. Since the metrics $\ell_{1}$ and
$\ell_{\infty}$ are dual, one hopes for interesting connections
between the two problems. Only very recently has this family of graphs
begun to be studied in \cite{VeomettRadcliffe}, in which the
vertex-isoperimetric problem is solved. 

\section{Monotonicity of the optimal boundary in isoperimetric problems on  ${\Z}^{k} \times {\N}^{d}$}

We denote vertices of a graph $G$ by $V(G)$, its edges by $E(G)$, and write $G=(V(G),E(G))$.  For $x \in V(G)$, we set $N_{V(G)}(x)=\{y \in V(G):\{x, y\} \in E(G) \}$ to be the vertex-neighbourhood of $x$ and $N_{E(G)}(x)=\{e \in E(G): e=\{x, y \} \in E(G) \text{ for some } y \in V(G) \}$ to be its edge neighbourhood. We write $N_G(x)$ for a neighbourhood of $x$ in $G$ when the distinction between vertex and edge does not matter.  \\

We will be needing the following definitions:

\begin{definition} \label{def central}
A graph $A=( {\Z}^{k+d},E(A))$ is locally symmetric if for every $x \in V(A)$, the neighbourhood $N_A(x)$ of $x$ is centrally symmetric about $x$. \\
A graph $G$ with $V(G)={\Z}^{k} \times {\N}^{d}$ is locally symmetric if there exists a graph $A=( {\Z}^{k+d},E(A))$ with centrally symmetric neighbourhoods such that for every $x \in V(G)$, the neighbourhood $N_G(x)$ of $x$ in $G$ is the intersection of the neighbourhood $N_A(x)$ of $x$ with $G$:  
\[
N_G(x)=N_A(x) \cap G.
\]
\end{definition}

\begin{definition}
A graph $G$  is induced by a $p$-norm, $1\leq p\leq \infty$, if there exists some constant $c$ such that for every $x \in V(g)$, $N_V(x)=\{y \in V(G): 0<\|x-y\|_p \leq c \}$.
\end{definition}

Note that if a graph $G$ with $V(G)={\Z}^{k} \times {\N}^{d}$ is induced by a $p$-norm, then it is homogenous.

Recall that colexicographical ordering is defined by 
\[
(a_1,a_2,\ldots,a_n) < (b_1,b_2,\ldots,b_n) \iff (\exists m>0)(\forall i > m) (a_i = b_i) \land (a_m < b_m)
\]

\begin{theorem}
\label{thm:The-minimum-is increasing in A} Let $G$ be a locally symmetric graph on ${\Z}^{k}$. The minimum edge-boundary is a monotonically increasing function of the volume.\end{theorem}

\proof
Let $B$ be a set of cardinality $|A|+1$ with optimal boundary.
We find a point which can be removed without increasing the boundary.
This implies that a set of cardinality $|A|$ has minimum boundary $|\partial A|$ less than
or equal to $|\partial B|$. Such a point has the property that it
has at least as many neighbours in ${\Z}^{k} \setminus B$ as
it does in $B$.
 We derive a contradiction by assuming that such a point does not
exist. 

Let
\[
R_{x}(y)=(2x_{1}-y_{1},2x_{2}-y_{2},\ldots,2x_{k}-y_{k})
\]
 denote the reflection of point $y$ in point $x$.

Let $H_m=\{x \in {\Z}^k: \|x\|_{\infty}=m\}$ and let $r=\max_r H_r \cap B \neq \emptyset$. Let $p \in H_r \cap B$ be greatest in colexicographical order. If $p$ has no neighbour in $B$, then it can clearly be removed without increasing the boundary. So let $y$ be a neighbour of $p$ in $B$. We show that the reflection $R_p(y)$ of $y$ in $p$ is not in $B$. Since $p$ is greater than $y$ colexicographically, $p_i=y_i$ for each $i>m$ and $p_m>y_m$ for some $m>0$. For each $i>m$, $2p_i-y_i=p_i$. Additionally,
\[
2p_m-y_m>p_m,
\]
so that $R_p(y)>p$ colexicographically. Therefore $R_p(y) \not \in B$. 
\proofend

\begin{theorem}
\label{thm:min perim} If $G$ is a graph on ${\Z}^{k} \times {\N}^{d}$ induced by a $p$-norm with constant $c<2$, then the minimum edge-boundary is a monotonically increasing function of the volume.\end{theorem}

\proof
Being induced by a $p$-norm, $G$ is locally symmetric. Therefore the argument from Theorem \ref{thm:The-minimum-is increasing in A} shows that for $b \in B$ of greatest colexicographical order, the reflection $R_b(y)$ of a neighbour $y \in B$ is outside of $B$. A new issue arises, that the image might land outside of the graph. For each index $i>k$ such that $2b_i-y_i <0$, we apply a reflection in the hyperplane $x_i=0$. The colexicographical order of the image can only increase, as its entries have increased. Moreover, each entry $2b_i-y_i <0$ maps to $y_i-2b_i>0$, so this image is outside of $B$ and inside the graph. Let $z$ be the image of $R_b(y)$ under these reflections. It remains to see that $z$ is a neighbour of $b$ and that this map is injective, i.e., no two neighbours $y$ and $y'$ map to the same point. 

To show that $z$ is a neighbour of $b$, it suffices to see that for each $i$, $|z_i-b_i|\leq |y_i-b_i|$. This inequality is clearly true for the indices $i$ for which no reflection in the axes occurs. So consider an index $i$ for which $2b_i-y_i <0$. We consider two cases, depending on whether $y_i < 3b_i$ or $y_i \geq 3 b_i$. \\ 
In the first case, $|z_i-b_i|=|y_i-3b_i|=3b_i-y_i$. On the other hand, $|y_i-b_i|=y_i-b_i$, so that
\[
|z_i-b_i| \leq |y_i-b_i| \iff y_i \geq 2b_i.
\] 
In the second case, $|z_i-b_i|=|y_i-3b_i|=y_i-3b_i$, so that
\[
|z_i-b_i| \leq |y_i-b_i| \iff b_i \geq 0.
\] 
Therefore $\|z-b\|_p \leq \|y-b\|_p$, so that $z$ is a neighbour of $b$. 

Next we show that if $y,y'$ are neighbours of $b$ then they do not map to the same point. Assume by contradiction that their images after reflections are the same. For each coordinate $i$, either
\[
2b_i-y_i=2b_i-y_i'
\]
or
\[
2b_i-y_i=y_i'-2b_i.
\]
Since $y \neq y'$, for at least one coordinate $i$, $y_i \neq y_i'$ and, consequently, $2b_i-y_i=y_i'-2b_i \implies 4b_i=y_i+y_i'$. Considering this equation over the nonnegative integers with constraint $y \neq y'$ shows that $|y_i'-b_i|\geq 2$ or $|y_i-b_i| \geq2$. It follows that one of $y$ and $y'$ is not a neighbour of $b$, a contradiction.
\proofend

Finally, we note that optimal sets are not necessarily nested. Consider the graph ${\N}^n$ induced by the $\infty$-norm with $c=1$. That is, the graph for which $x \sim y$  iff $\max_{i}|x_{i}-y_{i}|\leq1$. The set of edges is $E(G)=\{\{x,y\}:x,y\in\mathbb{N}^{n}\mbox{ and }\|x-y\|_{\infty}=\max_{i}|x_{i}-y_{i}|=1\}$. We consider the edge-isoperimetric problem on this
graph.

\begin{proposition}
The optimal sets of $\mathbb{N}^{2}$ are not nested.\end{proposition}

\begin{figure}[hbtp]
\centering
\includegraphics[scale=0.3]{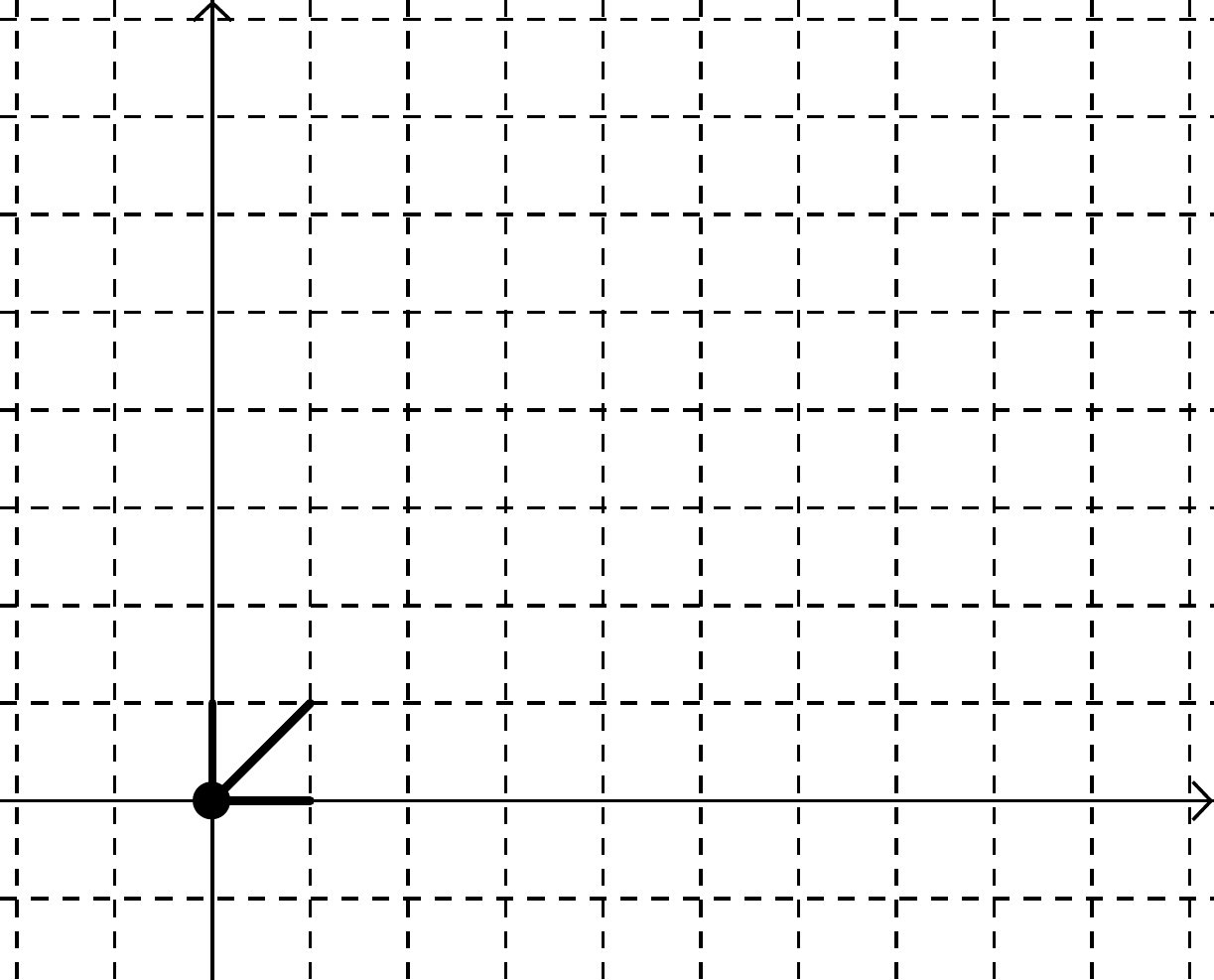} \includegraphics[scale=0.3]{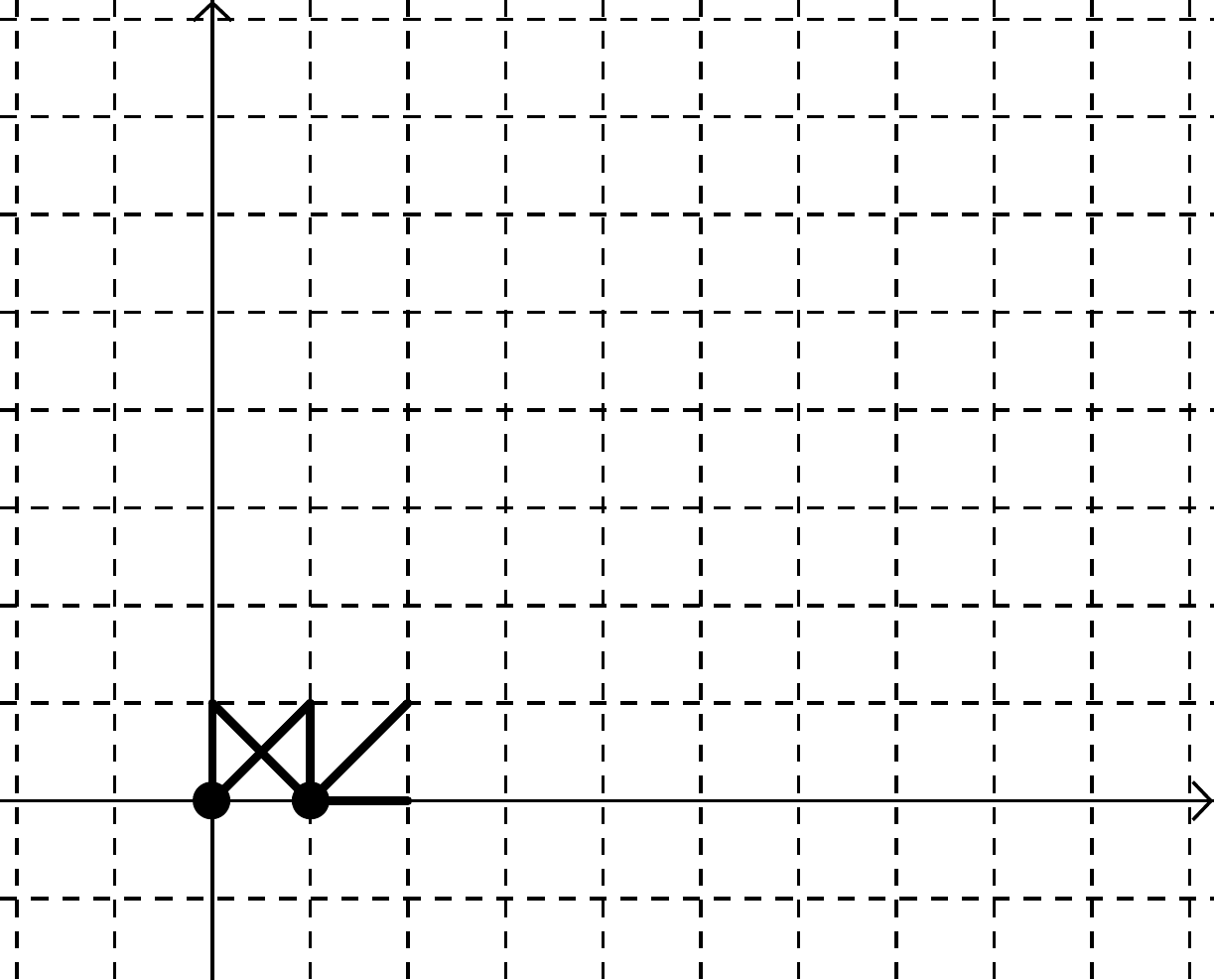}
\includegraphics[scale=0.3]{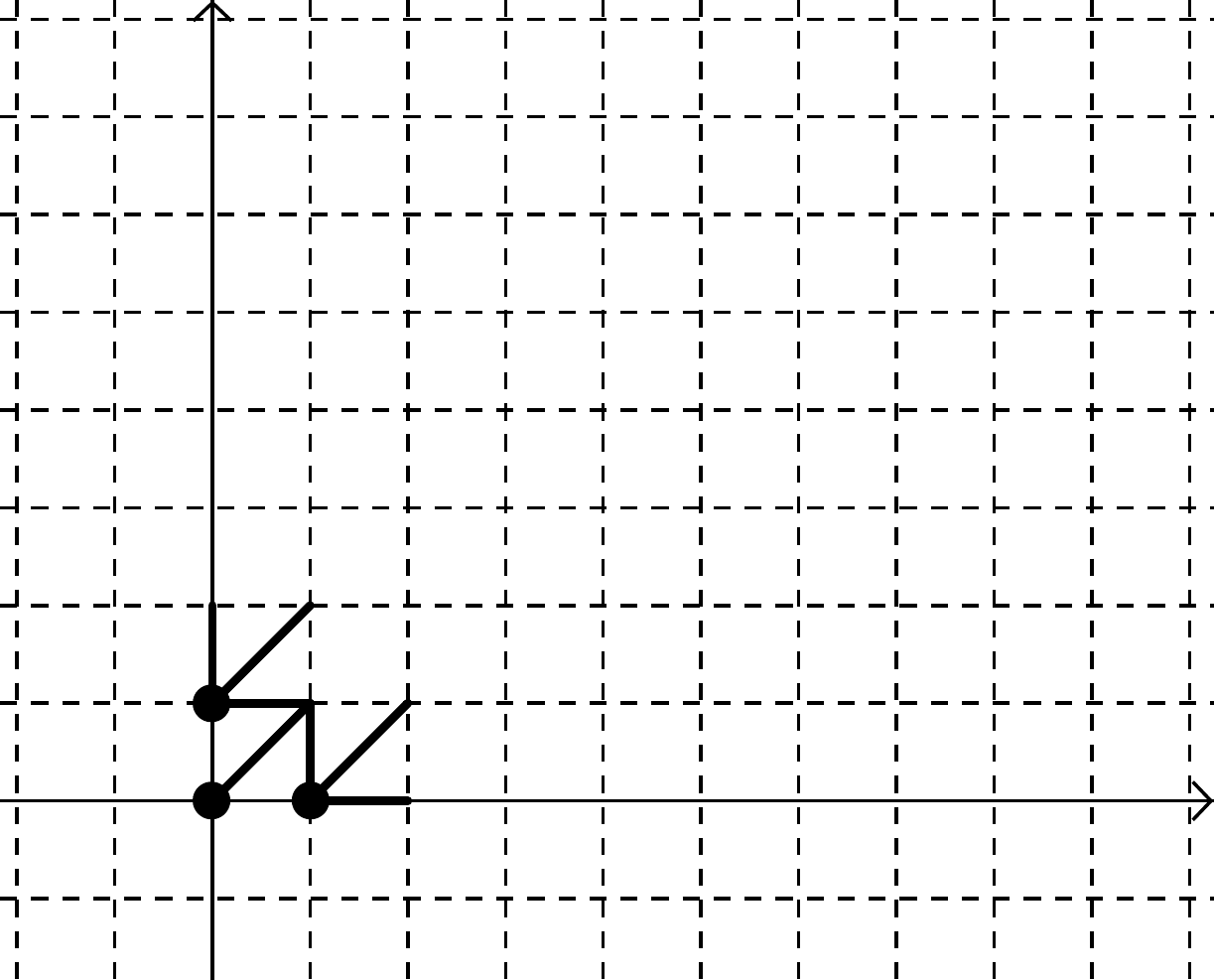}

\includegraphics[scale=0.3]{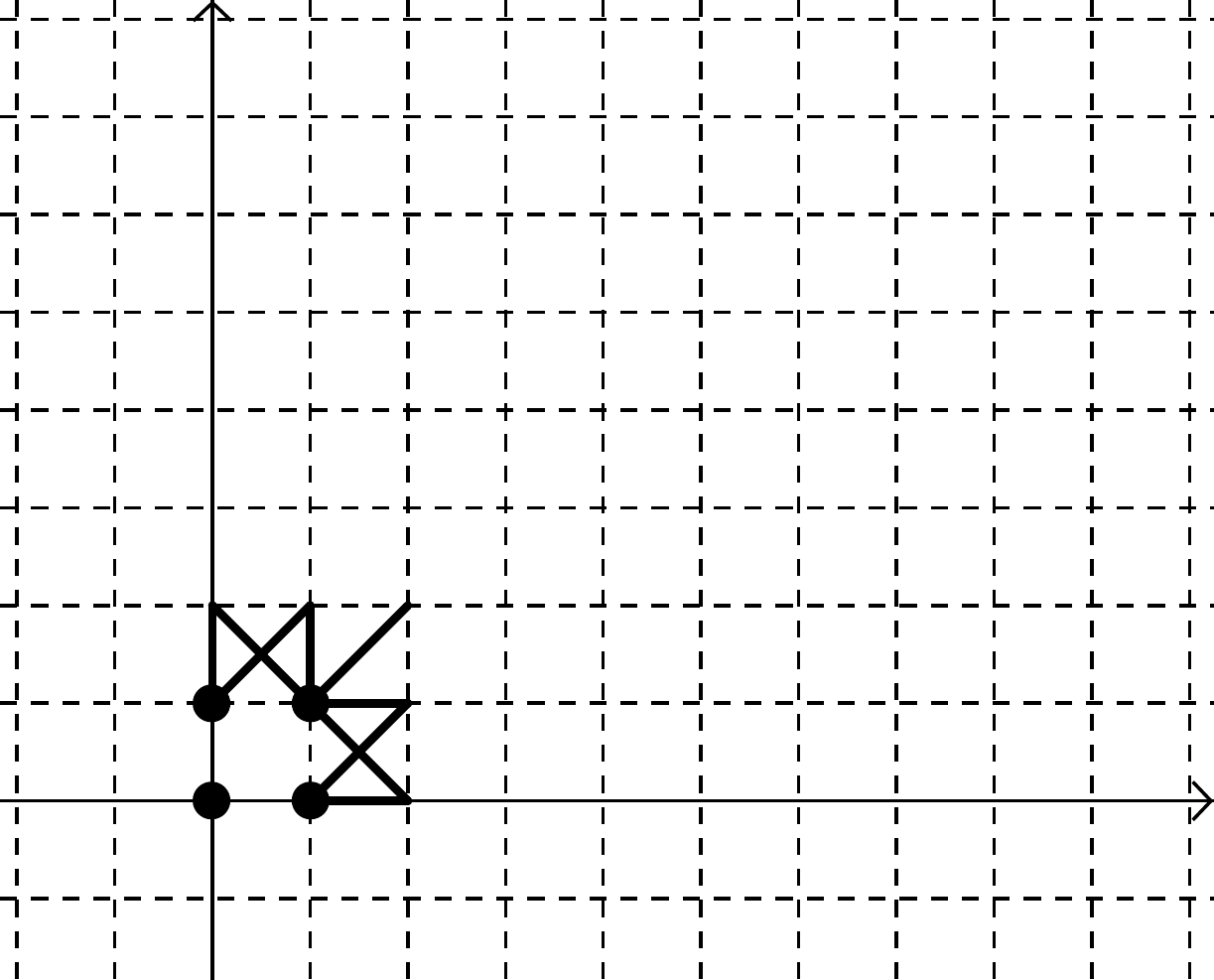} \includegraphics[scale=0.3]{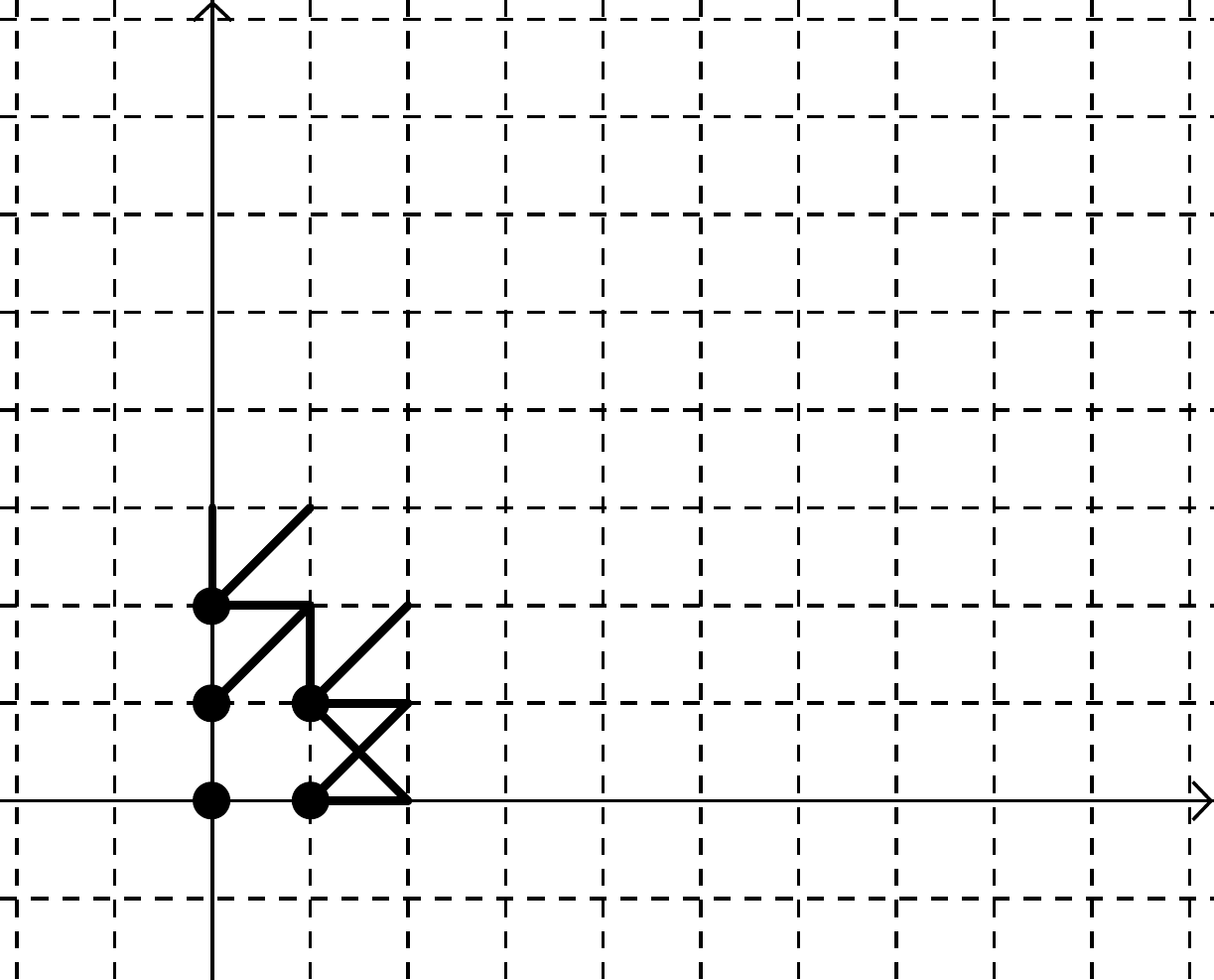}
\includegraphics[scale=0.3]{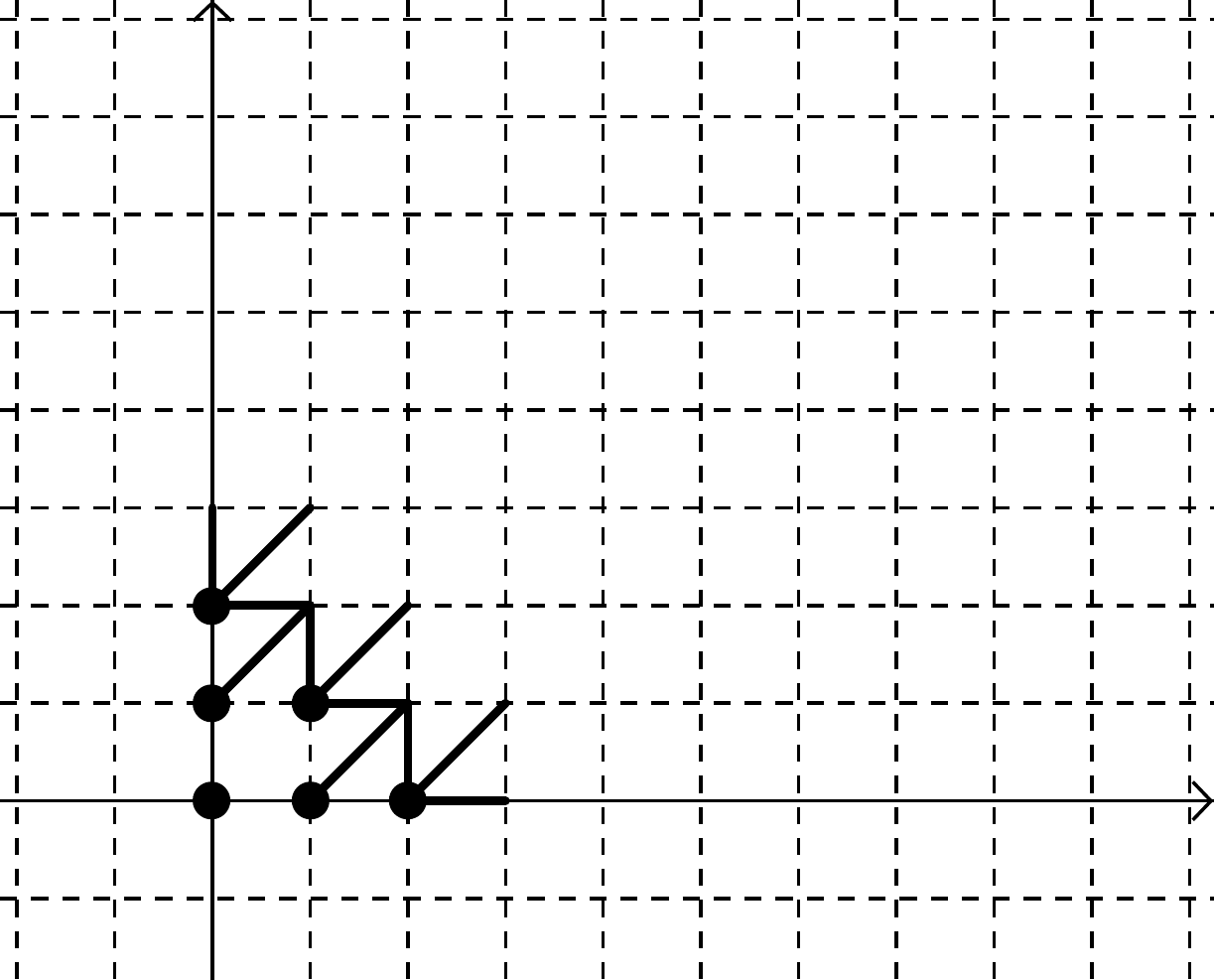}

\includegraphics[scale=0.3]{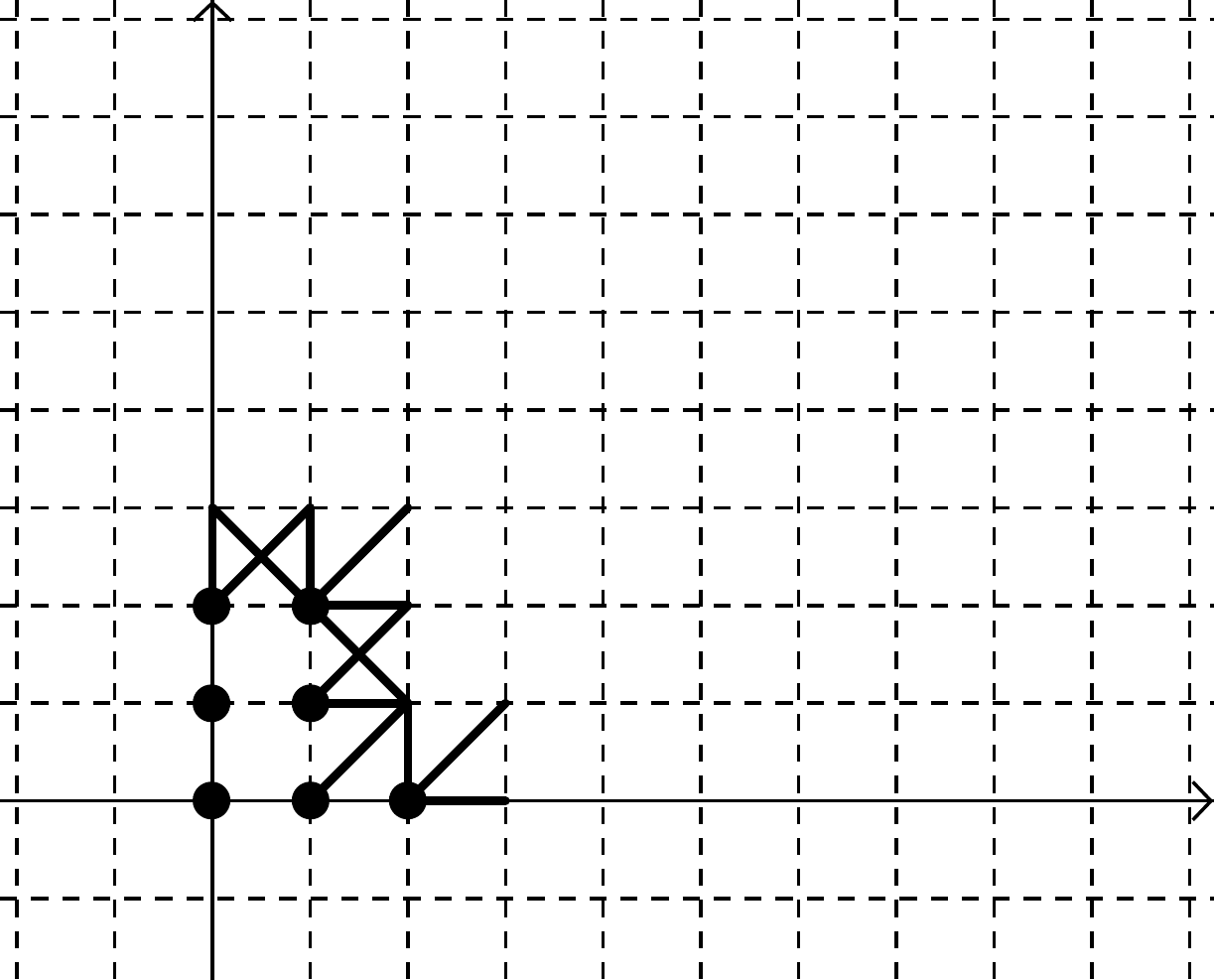} \includegraphics[scale=0.3]{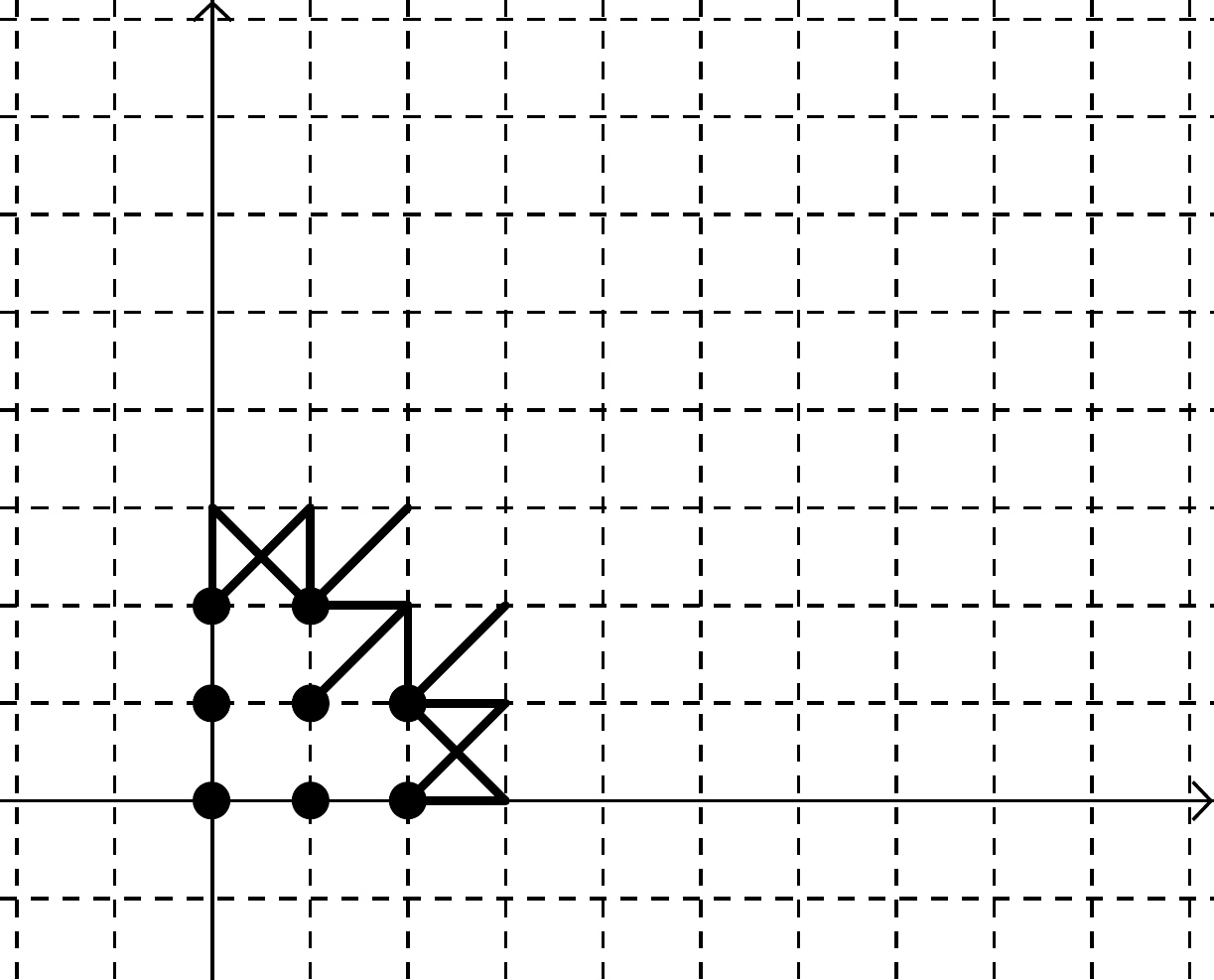}
\includegraphics[scale=0.3]{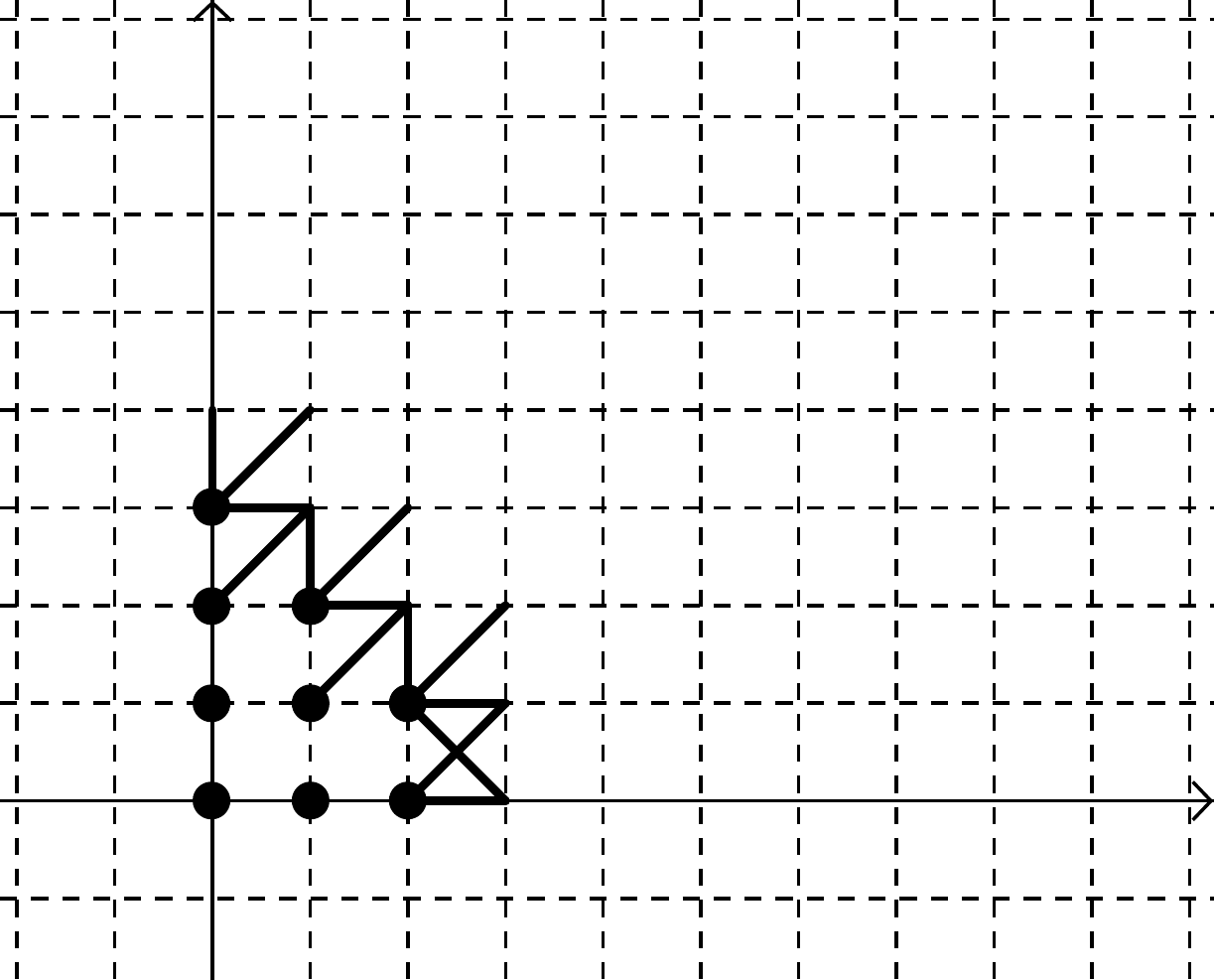}

\includegraphics[scale=0.3]{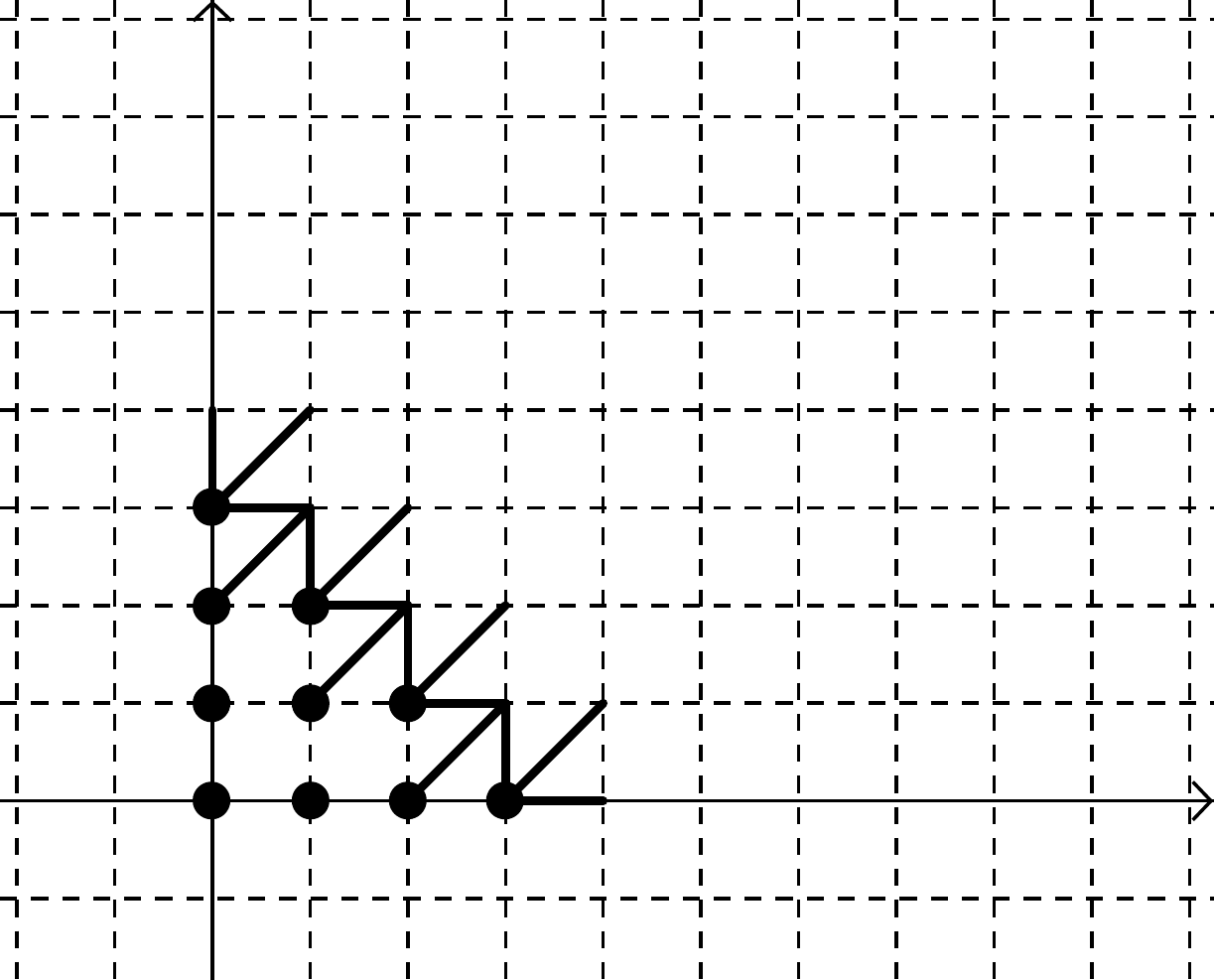} \includegraphics[scale=0.3]{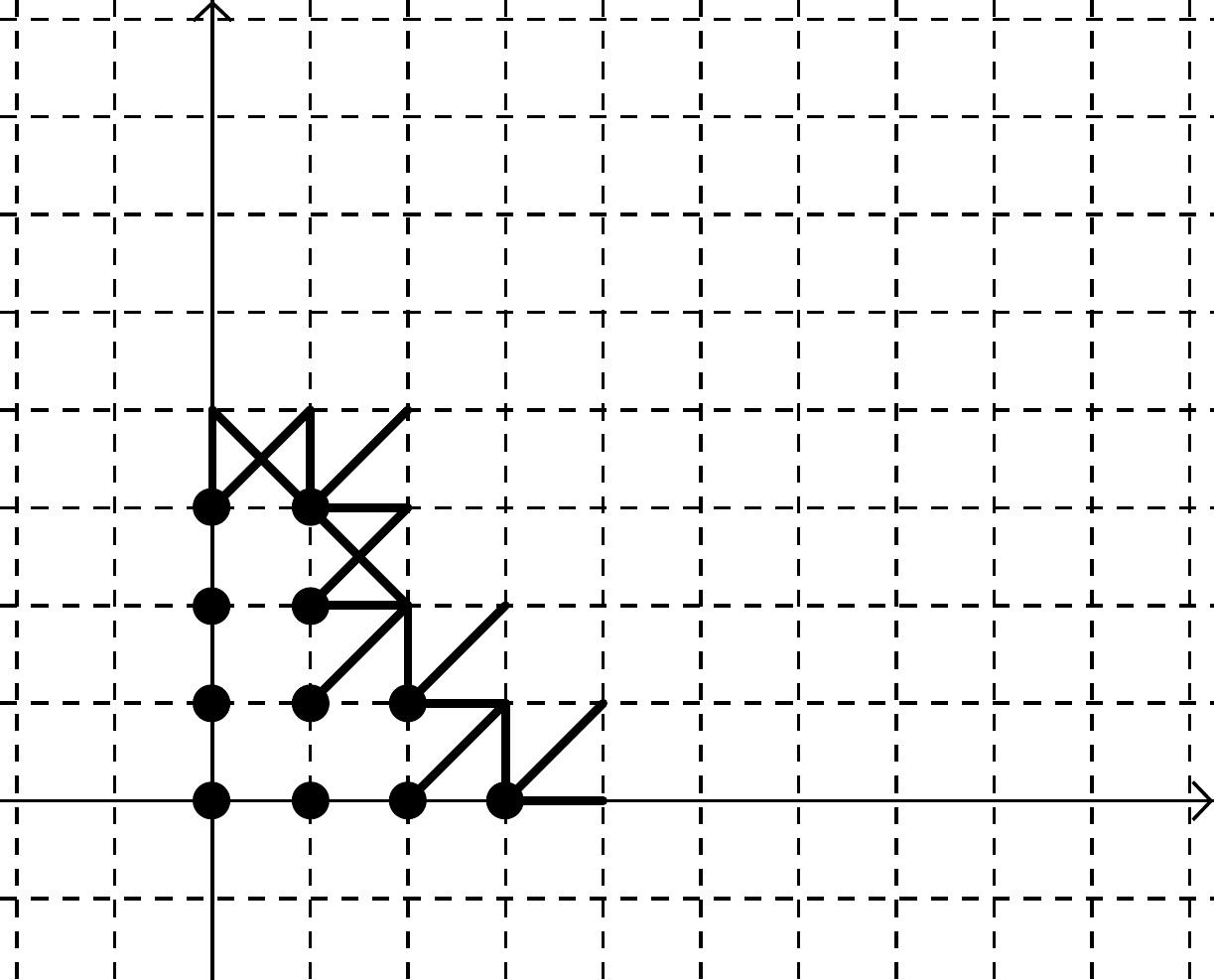}

\caption{\label{fig: nested sets}The unique (up to reflection in the line
$y=x$) sequence of optimal nested sets for $1\leq|A|\leq11$. The
optimal nested set with $|A|=11$ has $|\partial A|=17$.}
\end{figure}

\begin{figure}[hbtp] 
\centering
\includegraphics[width=2.4in]{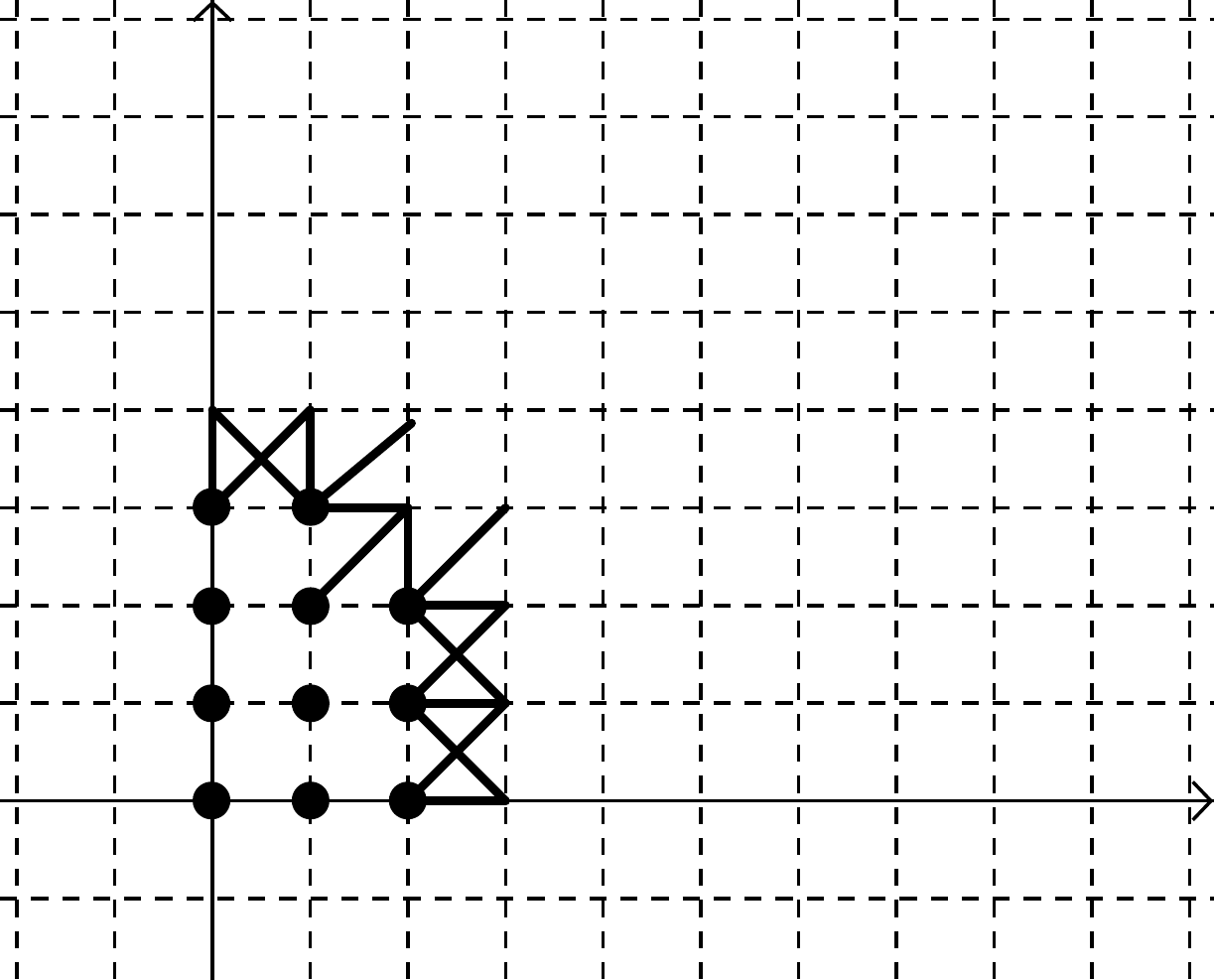}

\caption{\label{fig: optimal set with a=00003D00003D11}An optimal set having
$|A|=11$ has $|\partial A|=16$. This set is better than the one of the same volume in Figure \ref{fig: nested sets} }
\end{figure}

\proof
Figure \ref{fig: nested sets} shows the uniquely determined up to
reflection in $y=x$ sequence of nested optimal sets of $\mathbb{N}^{2}$.
Figure \ref{fig: optimal set with a=00003D00003D11} shows a set with
$|A|=11$ which has a smaller boundary than the optimal nested set
with $|A|=11$.
\proofend

Nested optimal sets are crucial for the technique of compression and
the fact that the optimal sets are not nested means that such an approach
will not be possible. However, as we will show, the monotonicity of the optimal boundary established here is a useful tool for obtaining bounds, proving optimality and understanding properties of optimal sets.

\section{The Edge-Isoperimetric Problem in $(\mathbb{N}^2,\infty)$}

Given a set in ${\N}^2$, it can be made connected without increasing its volume by translating the connected components towards the origin. Moreover, it can be made to touch both axes. Call $A$ the resulting set. Let $X$ be the point on the $x$-axis with greatest $x$-coordinate and let $Y$ be the point on the $y$-axis of greatest $y$-coordinate. There is a connected subset $C$ of $A$ containing $X$ and $Y$. We argue that all points bounded by $C$ and the axes are in $A$ if $A$ is optimal.

\begin{definition} Given a connected subset $C$ of ${\N}^2$ containing at least one point on the $x$-axis and one point on the $y$-axis, we say that a point $z$ is bounded by $C$ if point $z$ is bounded by some piecewise linear curve formed by a subset of the edges of $C$ and the axes.
\end{definition}  

\begin{lemma} \label{fillingCurves}
If $A$ is an optimal set, then it contains all lattice points bounded by the axes and a connected subset $C$ containing points $X$ and $Y$.
\end{lemma}

\begin{figure}[hbtp] 
\centering
\includegraphics[width=2.4in]{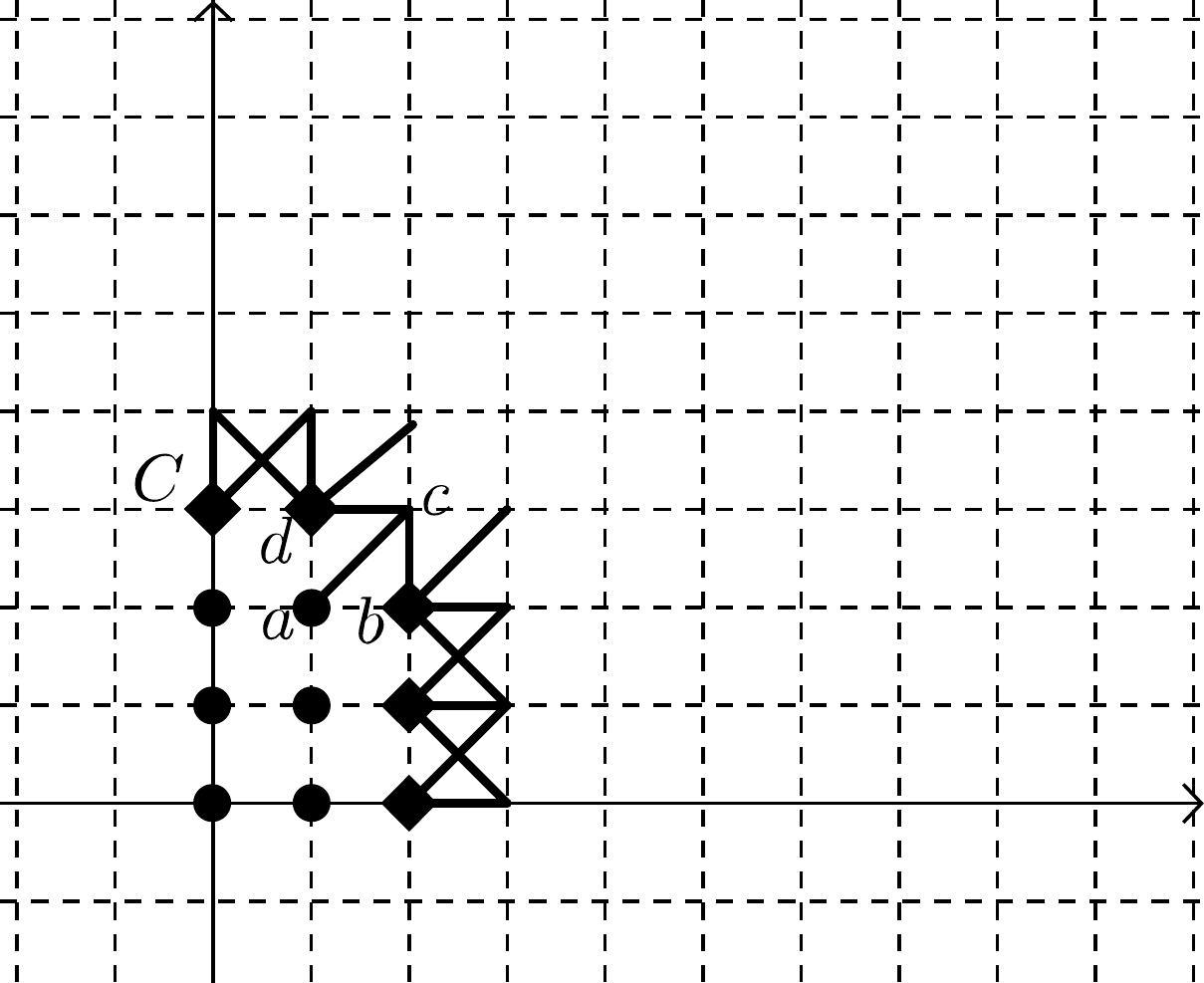}

\caption{\label{fig:outer layer} A connected subset $C$ containing $X$ and $Y$ is indicated by diamonds. The perimeter of the set bounded by $C$ receives a contribution not only from the points of $C$, but also from the layer adjacent to $C$, e.g., point $a$. }
\end{figure}

\proof
Assume otherwise. Consider the subset $B_C \subset {\N}^2$ of points lying on or below $C$ and let $B=B_C \cup A$. Let $A_C=B_C \cap A$ be the points of $A$ lying on or below $C$.  By assumption, $|B_C| >|A_C|$, so that $|B|>|A|$. We show that $|\partial B_C|<|\partial A_C|$. This implies that $|\partial B|<|\partial A|$, contradicting Theorem \ref{thm:min perim}. \\
The perimeter $\partial B_C$ of $B_C$ consists not only of all edges of $C$ which lie outside of $B_C$, but also of the edges of the outermost layer of points of $B_C \setminus C$ (see Figure \ref{fig:outer layer}) if such points exist. In particular, a point $a \in B_C \setminus C$ contributes whenever the following conditions are met. Point $a$ is the vertex of a unit square $abcd$, the opposite corner $c$ is in the complement of $B_C$, and the other two corners $b$ and $d$ are on $C$. In this case, edge $ac$ is added to the perimeter. Suppose that such a point $a \in B_C \setminus A_C$ exists. Then the addition of $a$ to $A_C$ adds a diagonal edge $ac$ to the count but takes away the two edges $ab$ and $ad$. If $a$ is part of more than one square, it is easy to see that the perimeter will still be strictly improved. If $B_C \setminus A_C$ contains some point which does not meet the conditions, then it does not contribute to the perimeter of $B_C$. Since there is at least one point of this form or of the prior form, filling in these points strictly improves the perimeter. This completes the proof. \proofend

\begin{conjecture} \label{setBelowCurve}
For every volume, there is an optimal set $A$ which consists of the points bounded by a connected subset $C \subset A$ touching the axes. 
\end{conjecture}

We call sets consisting of the points bounded by a connected subset which touches the axes \emph{bounded}. We will now investigate bounded sets. Though there might be sets which are better, we will still be able to learn much about the problem by considering bounded sets.

\begin{definition}
Point $g \in {\Z}^n$ is a $j$-gap of a set $A$ if $g\notin A$ and there exists some point $p \in A$ with $p_{i}=g_{i}\mbox{ }\forall i\neq j$
and $g_{j}<p_{j}$. 
\end{definition}

We will say that a set has no gaps with it has no $j$-gaps for any $j=1,2,\ldots,n$.

\begin{figure}[hbtp] 
\centering
\includegraphics[width=2.4in]{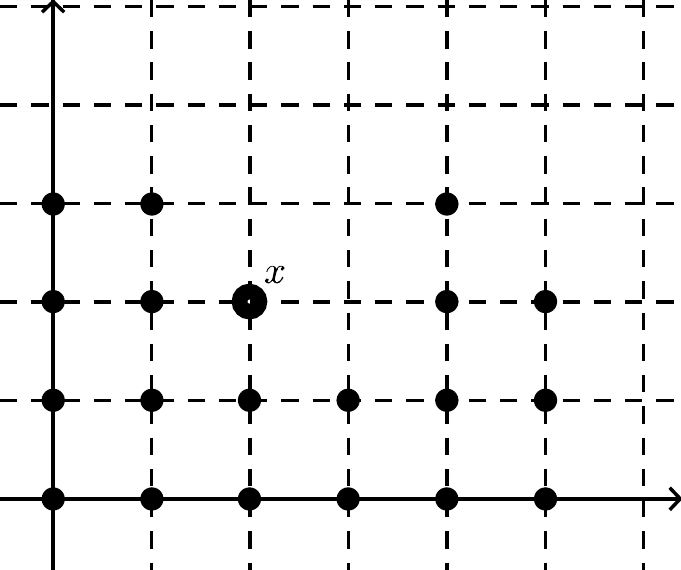}

\caption{\label{fig:convex} A bounded optimal set $A$ cannot have gaps since filling in the gap with points such as $x$ decreases the boundary. }
\end{figure}

\begin{theorem} \label{convexOptimalSet}
For every volume, a bounded set can be modified to a bounded set with no gaps without increasing the boundary. 
\end{theorem}
\proof
We fill in the $1$-gaps starting from the lowest gaps by adding points such as $x$ in Figure \ref{fig:convex}, all the while decreasing the perimeter. By Theorem \ref{thm:min perim}, the optimal perimeter is a monotonic increasing function of the volume, so we see that there can be no $1$-gaps. The same argument applies to the $2$-gaps. \proofend

Let $A_t=\{x \in A: x_1=t\}$.

\begin{lemma}
\label{lem: columns heights decrease by less than 2}A bounded optimal set
can be chosen to have at most one $t\in\{1,2,...,k-1\}$ for which
$|A_{t}|-|A_{t+1}|\geq2$, and this $t$ can be chosen to be $k-1$.\end{lemma}
\proof
Let $t_{0}$ be the first $t$ for which $|A_{t}|-|A_{t+1}|\geq2$
and assume that $t_{0}\neq k-1$. We can then transfer points from
$A_{k}$ to $A_{t_{0}+1}$ until $|A_{t_{0}}|-|A_{t_{0}+1}|=1$ or
$A_{k}$ has no more points left. Throughout, the perimeter does not
increase because any point transferred shared at most $8$ edges with
other points, and after the transfer shares at least $8$ edges. If
$|A_{t_{0}}|-|A_{t_{0}+1}|$ is still greater than $2$, then we can
transfer points from column $A_{k-1}$, and continue in this way until
either $|A_{t_{0}}|-|A_{t_{0}+1}|=1$ or the number of columns has
been reduced so that $t_{0}+1$ is the last column. 
\proofend

Let us now view the heights of the columns as a function $h:\mathbb{N}\rightarrow\mathbb{N}$
given by $h(t)=|A_{t}|$. For brevity we will say that $h$ is constant
whenever we mean that $h$ is constant on its support. We have already
shown in Theorem \ref{convexOptimalSet} that a bounded optimal
set exists which has $h$ non-increasing. We will show that $h$ can
be made to take on a specific form. Before we do that, however, we
note some special cases that are the only exceptions to the following
Lemma. If $|A|=1$, $|A|=2$ or $|A|=4$, then it is easy to see that
for the optimal set $h$ is constant. These are the only cases in
which $h$ will be constant, as we show next.
\begin{lemma}
\label{lem: shape of optimal sets}Without increasing the boundary,
a bounded set $A$ can be transformed into a bounded set $B$ for which $h_{B}$ is constant
on $\{1,2,...,c-1\}$ and strictly decreasing on $\{c,c+1,...,k\}$.
Moreover, if $|A|\neq1,2,4$, we can choose $c<k$.\end{lemma}
\proof
Let $S=\{x,x+1,...,x+l\}$ be a maximal set of least $x$ for which
$h(x+i)<h(x+i-1)\mbox{ }\forall i=1,...,l$. Assume further that $x+l<k$.
We consider two cases: when $x+l=k-1$ and when $x+l<k-1$. In the
former case, we take the top point from column $A_{k}$, which has
at most $6$ shared edges, and place it on top of column $A_{k-1}$,
and now it has at least $6$ shared edges. This reduces us to the
second case. In this case, we take a point from column $A_{k}$, which
can have at most $8$ shared edges, and place it at the top of column
$A_{x+l}$. Because $h(x+l)=h(x+l+1)$ and $h(x+l-1)=h(x+l)+1$, the
new point has $8$ shared edges, so the perimeter is not increased.
We have now reduced $S$ to $S\setminus\{x+l\}$, and we continue
inductively. This shows that $A$ can be transformed into a set $B$
for which $h_{B}$ is constant on $\{1,2,...,c-1\}$ and strictly
decreasing on $\{c,c+1,...,k\}$.

To see that we can choose $c<k$, assume that $|A|\neq1,2,4$ and
that $h$ is constant. If $k=1$, then we can take the top point of
$A_{k}$ and place it at $A_{2}$ without increasing the perimeter.
If $k=2$, we can take the point $P$ at the top of $A_{k}$ and place
it at $A_{k+1}$. Since $h_{A}$ was constant, $P$ had at most $3$
neighbours. Placing $P$ at $A_{k+1}$ guarantees $2$ neighbours
and $3$ boundary edges, so the perimeter is not increased. For $k\geq3$,
we can take the point $P$ from the top of $A_{k}$ and place it at
the top of $A_{1}$, and because $k\geq3$, $P$ will not have $2$
neighbours and $3$ boundary edges.
\proofend

\begin{figure}[hbtp]
\centering
\includegraphics[width=2.6in]{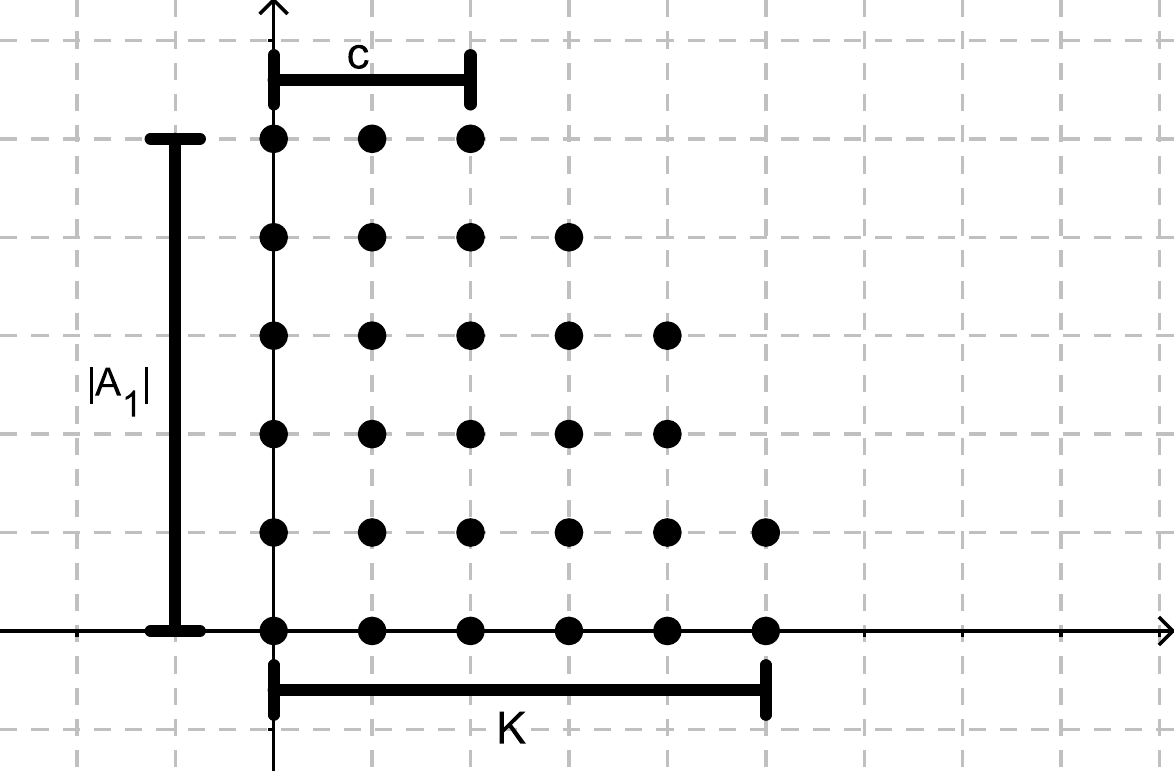}
\caption{Schematic depicting the form of an optimal bounded set guaranteed by Lemmas
\ref{lem: columns heights decrease by less than 2} and \ref{lem: shape of optimal sets}.
The first $c$ columns have the same height $|A_{1}|$, the next $k-c-1$
columns have heights decreasing by $1$ at each step, and the last
column $A_{k}$ has height $|A_{k}|$ less than $|A_{k-1}|$.}
\end{figure}

We can now find the perimeter of a general set subject to the conditions
in the Lemmas above in terms of $|A_{1}|$, $|A_{k}|$, $k$ and $c$.
There are $|A_{1}|$ horizontal edges and $k$ vertical edges. There
are $\sum_{t=1}^{k-1}(|A_{t}|-|A_{t+1}|+1)+|A_{k}|=|A_{1}|+k-1$ edges
parallel to $e_{1}+e_{2}$, $|A_{k}|-1+\max\{|A_{k-1}|-|A_{k}|-1,0\}=|A_{k}|-1+|A_{k-1}|-|A_{k}|-1=|A_{k-1}|-2$
edges in the $e_{1}-e_{2}$ direction, and $\sum_{t=2}^{k}\delta_{|A_{t}|,|A_{t-1}|}=c-1$
edges in the direction $e_{2}-e_{1}$. Consequently, the perimeter
is $2|A_{1}|+|A_{k-1}|+c+2k-4=2|A_{1}|+|A_{1}|-(k-1-c)+c+2k-4$, which
is equal to 
\begin{equation}
|\partial A|=3|A_{1}|+2c+k-3.\label{eq:perimeter}
\end{equation}

We also know that $\sum_{t=1}^{k}|A_{t}|=|A|$. Therefore $|A|=c|A_{1}|+\sum_{i=1}^{k-c-1}(|A_{1}|-i)+|A_{k}|,$
which simplifies to 
\begin{equation}
|A|=(k-1)|A_{1}|+|A_{k}|-\frac{(k-c-1)(k-c)}{2}.\label{eq:volume}
\end{equation}

Combining Theorem \ref{thm:min perim} on the
monotonicity of the perimeter and equation \ref{eq:perimeter}, we
obtain:
\begin{corollary}
Let $A$ be a bounded optimal set with $|A_{k}|<|A_{k-1}|-1$. Then the optimal
perimeter of bounded sets of cardinality $|A|+1$ is $|\partial A|$ and a bounded
optimal set of cardinality $|A|+1$ is given by $A\cup\{(k,|A_{k}|+1)\}$.
\end{corollary}

\begin{example}
We show that simplices are not always optimal.
Consider the simplex given by $h_{i}=15-i$, $i=1,2,...,14$. By increasing
$|A_{1}|$ by $1$ while preserving the shape from Lemma \ref{lem: shape of optimal sets},
we get a truncated simplex given by $h_{i}=16-i$, $i=1,2,...,10$.
The perimeter of the simplex is $56$, whereas that of the truncated
simplex is $55$.
\end{example}

Given $|A_{1}|$, $c$ and $|A|$, the set $A$ is determined. Indeed,
we know that 
\[
|A_{i}|=\begin{cases}
|A_{1}| & i\leq c\\
|A_{1}|-i+c & c+1\leq i\leq k-1\\
|A|+\frac{(k-c-1)(k-c)}{2}-(k-1)|A_{1}| & i=k
\end{cases}
\]

Moreover, $k$ is the unique positive integer such that 
\[
\sum_{i=1}^{k-c-1}(|A_{1}|-i)<|A|-c|A_{1}|\leq\sum_{i=1}^{k-c}(|A_{1}|-i).
\]

Solving for $k$, we obtain 
\[
k=\lceil\frac{1}{2}(-1+2|A_{1}|-\sqrt{1+8(\binom{|A_{1}|}{2}-|A|+c|A_{1}|)}\rceil+c.
\]

Therefore the problem is to minimize 
\[
|\partial A|=3|A_{1}|+3c+\lceil\frac{1}{2}(-1+2|A_{1}|-\sqrt{1+8(\binom{|A_{1}|}{2}-|A|+c|A_{1}|)})\rceil-3
\]
 
\[
=4|A_{1}|+3c-3-\lfloor\frac{1}{2}(1+\sqrt{1+8(\binom{|A_{1}|}{2}-|A|+c|A_{1}|)})\rfloor.
\]

\begin{lemma}
Any bounded set $A$ can be transformed into one for which $|A_{1}|\geq c$ without increasing the boundary.\end{lemma}
\proof
Assume that $|A_{1}|<c$. We reflect $A$ in the line $y=x$ to obtain
a new set $B$ which has $|A_{1}|$ columns. The first $|A_{k}|$
are of height $k$. Columns $|A_{k}|+1$ through $|A_{k-1}|$ are
of height $k-1$. The remaining columns decreased in height by steps
of $1$, with column $|A_{k-1}|+1$ having height $k-2$ and column
$|A_{1}|$ having height $c$. Since $c>|A_{1}|>|A_{k-1}|-|A_{k}|$,
we can take points from column $|A_{1}|$ of height $c$ and place
them on top of columns $|A_{k}|+1$,...,$|A_{k-1}|-1$ without increasing
the perimeter. The resulting set has the form of Lemma \ref{lem: shape of optimal sets}.
The new parameters are $\tilde{|A_{1}|}=k$, $\tilde{k}=|A_{1}|$
and $\tilde{c}=|A_{k-1}|-1=|A_{1}|-k+c$. Then 
\[
c>|A_{1}|\implies\tilde{c}=|A_{1}|-k+c<2c-k.
\]

Since $k>c$, $|\tilde{A_{1}}|>\tilde{c}.$ 
\proofend

In order to obtain a lower bound for the sets considered, we relax our problem to a continuous
one:

\begin{equation*}
\begin{aligned}
& \underset{|A_{1}|,c\in\mathbb{R}}{\text{minimize}}
& & 4|A_{1}|+3c-3-\frac{1}{2}(1+\sqrt{1+8(\binom{|A_{1}|}{2}-|A|+c|A_{1}|)}) \\
& \text{subject to} 
& & 1\leq|A_{1}|\leq|A|, \\
&&& \frac{|A|-\binom{|A_{1}|}{2}}{|A_{1}|}\leq c\leq|A_{1}|.\;
\end{aligned}
\end{equation*}

For any $|A|\geq2$, we can establish via a direct calculation that
the minimum value of the objective is $\sqrt{\frac{7}{2}}\sqrt{8|A|-1}-2$
given by the unconstrained minimizer $|A_{1}|=\frac{3\sqrt{8|A|-1}}{2\sqrt{14}}$
and $c=\frac{1}{28}(14+\sqrt{14}\sqrt{8|A|-1})$, and the value is better
than the value of the function on the boundary of the feasible region.

To obtain an upper bound, we will utilize the monotonicity of the
perimeter from Theorem \ref{thm:min perim}.
Let $m\in\mathbb{N}$ and set $|A|^{*}=7m^{2}$, $|A_{1}|^{*}=3m$
and $c^{*}=m$. It is again a simple calculation to verify that these
values give a feasible point. The function 
\[
g(|A_{1}|,c)=4|A_{1}|+3c-2-\frac{1}{2}(1+\sqrt{1+8(\binom{|A_{1}|}{2}-|A|+c|A_{1}|)})
\]

is an upper bound for the perimeter. For $|A|^{*}=7m^{2}$, 
\[
g(|A_{1}|^{*},c^{*})=15m-\frac{1}{2}\sqrt{4m^{2}-12m+1}-\frac{5}{2}
\]
 
\[
=\frac{15}{\sqrt{7}}\sqrt{|A|^{*}}-\frac{1}{2}\sqrt{\frac{4}{7}|A|^{*}-\frac{12}{\sqrt{7}}\sqrt{|A|^{*}}+1}.
\]

For a general $|A|$, we find $m$ such that $7(m-1)^{2}<|A|\leq7m^{2}$.
Then $|A|^{*}\leq|A|+2\sqrt{7|A|}-8$, so that 
\[
|\partial A|\leq\frac{15}{\sqrt{7}}\sqrt{|A|+2\sqrt{7|A|}-8}-\frac{1}{2}\sqrt{\frac{4}{7}(|A|+2\sqrt{7|A|}-8)-\frac{12}{\sqrt{7}}\sqrt{|A|+2\sqrt{7|A|}-8}+1}.
\]

This complicated expression is asymptotically 
\[
\frac{15}{\sqrt{7}}\sqrt{|A|}-\frac{1}{\sqrt{7}}\sqrt{|A|}=\sqrt{\frac{7}{2}}\sqrt{8|A|}.
\]

Note, however, that the upper bound is real if and only if $|A|\geq36$. 
\begin{theorem}
\label{thm:Lower and Upper Bounds}Let the cardinality of a bounded optimal set $A$
be $|A|\geq36$. Then the perimeter $|\partial A|$ is bounded below by
\[
\lceil\sqrt{\frac{7}{2}}\sqrt{8|A|-1}-2\rceil
\]

and above by

\[ \lfloor\frac{15}{\sqrt{7}}\sqrt{|A|+2\sqrt{7|A|}-8}-\frac{1}{2}\sqrt{\frac{4}{7}(|A|+2\sqrt{7|A|}-8)-\frac{12}{\sqrt{7}}\sqrt{|A|+2\sqrt{7|A|}-8}+1}\rfloor.
\]

Moreover, the difference between the upper and lower bound does not
exceed the constant $\frac{35}{2}$.\end{theorem}
\proof
The upper and lower bounds were shown above. Rather than consider
the upper bound as it stands, we consider the slightly weaker but
simpler upper bound $u(|A|)$ equal to
\[ \frac{15}{\sqrt{7}}\sqrt{|A|+2\sqrt{7|A|}-8}-\frac{1}{2}\sqrt{\frac{4}{7}(|A|+2\sqrt{7|A|}-8)-\frac{12}{\sqrt{7}}\sqrt{|A|+2\sqrt{7|A|}-8}}
\]
obtained by dropping the $1$ inside of the square root. Then a calculation
shows that the difference $d(|A|)$ between the upper bound $u(|A|)$
and lower bound $l(|A|)=\sqrt{\frac{7}{2}}\sqrt{8|A|-1}-2$ has a
non-vanishing derivative. Moreover, at $|A|=39$, the first point
at which this upper bound is defined, the derivative of the difference
is positive, so that the difference is an increasing function. Taking
the limit, we obtain the value $\frac{35}{2}$. For any $36\leq|A|\leq38$,
a direct calculation shows that the difference is at most $\frac{35}{2}$.
\proofend

Note that the growth of the boundary, even for bounded optimal sets, is slower than linear, though by Theorem \ref{thm:min perim}
the perimeter is an increasing function. Therefore 
\begin{corollary}
There exist arbitrarily long consecutive values of the volume for
which the minimum boundary is the same. 
\end{corollary}

\section{Acknowledgements}

I would like to thank Professor Simon Brendle for his helpful suggestions throughout the writing of this note.

\end{document}